\documentclass[conference]{IEEEtran}
\IEEEoverridecommandlockouts
\usepackage{cite}
\usepackage{subcaption}
\usepackage{amsmath,amssymb,amsfonts}
\usepackage{algorithmic}
\usepackage{graphicx}
\usepackage{textcomp}
\usepackage{url}
\usepackage{xcolor}
\def\BibTeX{{\rm B\kern-.05em{\sc i\kern-.025em b}\kern-.08em
    T\kern-.1667em\lower.7ex\hbox{E}\kern-.125emX}}
\begin{document}


\newenvironment{Abstract}
  {\begin{center}
    \textbf{Abstract} \\
    \medskip
  \begin{minipage}{4.8in}}
  {\end{minipage}
\end{center}}

\def\noproof{\hspace{1em}\blackslug}
\def\mm#1{\hbox{$#1$}}                 
\def\mtx#1#2{\renewcommand{\arraystretch}{1.2}%
    \left(\! \begin{array}{#1}#2\end{array}\! \right)}
			\newcommand{\mQ}{\mathcal{Q}}
\newcommand{\ie}{i.e.\@\xspace} 
\newcommand{\eg}{e.g.\@\xspace} 
\newcommand{\etal}{et al.\@\xspace} 
\newcommand{\Matlab}{\textsc{Matlab}\xspace}
\newcommand{\e}[2]{{\small e}$\scriptstyle#1$#2}
\renewcommand{\vec}[1]{\ensuremath{\mathbf{#1}}}  
\newcommand{\abs}[1]{\ensuremath{|#1|}}
\newcommand{\Real}{\ensuremath{\mathbb{R}}}
\newcommand{\inv}{^{-1}}
\newcommand{\dom}{\mathrm{dom}}
\newcommand{\minimize}[1]{\displaystyle\minim_{#1}}
\newcommand{\minim}{\mathop{\hbox{\rm minimize}}}
\newcommand{\maximize}[1]{\displaystyle\maxim_{#1}}
\newcommand{\maxim}{\mathop{\hbox{\rm maximize}}}
\newcommand{\twonorm}[1]{\norm{#1}_2}
\def\onenorm#1{\norm{#1}_1}
\def\infnorm#1{\norm{#1}_{\infty}}
\def\disp{\displaystyle}
\def\m{\phantom-}
\def\EM#1{{\em#1\/}}
\def\ur{\mbox{\bf u}}
\def\FR{I\!\!F}
\def\Complex{I\mskip -11.3mu C}
\def\sp{\hspace{10pt}}
\def\subject{\hbox{\rm subject to}}
\def\eval#1#2{\left.#2\right|_{\alpha=#1}}

\def\tmat#1#2{(\, #1 \ \ #2 \,)}
\def\tmatt#1#2#3{(\, #1 \ \  #2 \ \  #3\,)}
\def\tmattt#1#2#3#4{(\, #1 \ \  #2 \ \  #3 \ \  #4\,)}

\def\mat#1#2{(\; #1 \quad #2 \;)}
\def\matt#1#2#3{(\; #1 \quad #2 \quad #3\;)}
\def\mattt#1#2#3#4{(\; #1 \quad #2 \quad #3 \quad #4\;)}
\def\blackslug{\hbox{\hskip 1pt \vrule width 4pt height 6pt depth 1.5pt
  \hskip 1pt}}
\def\boxit#1{\vbox{\hrule\hbox{\vrule\hskip 3pt
  \vbox{\vskip 3pt #1 \vskip 3pt}\hskip 3pt\vrule}\hrule}}
\def\cross{\scriptscriptstyle\times}
\def\cond{\hbox{\rm cond}}
\def\det{\mathop{\hbox{\rm det}}}
\def\diag{\mathop{\hbox{\rm diag}}}
\def\dim{\mathop{\hbox{\rm dim}}}  
\def\exp{\mathop{\hbox{\rm exp}}}  
\def\In{\mathop{\hbox{\it In}\,}}
\def\boundary{\mathop{\hbox{\rm bnd}}}
\def\interior{\mathop{\hbox{\rm int}}}
\def\Null{\mathop{\hbox{\rm null}}}
\def\rank{\mathop{\hbox{\rm rank}}}
\def\op{\mathop{\hbox{\rm op}}}
\def\Range{\mathop{\hbox{\rm range}}}
\def\range{\mathop{\hbox{\rm range}}}
\def\sign{\mathop{\hbox{\rm sign}}}
\def\sgn{\mathop{\hbox{\rm sgn}}}
\def\Span{\mbox{\rm span}}
\def\trace{\mathop{\hbox{\rm trace}}}
\def\drop{^{\null}}
\def\float{fl}
\def\grad{\nabla}
\def\Hess{\nabla^2}
\def\half  {{\textstyle{1\over 2}}}
\def\third {{\textstyle{1\over 3}}}
\def\fourth{{\textstyle{1\over 4}}}
\def\sixth{{\textstyle{1\over 6}}}
\def\Nth{{\textstyle{1\over N}}}
\def\nth{{\textstyle{1\over n}}}
\def\mth{{\textstyle{1\over m}}}
\def\invsq{^{-2}}
\def\limk{\lim_{k\to\infty}}
\def\limN{\lim_{N\to\infty}}
\def\mystrut{\vrule height10.5pt depth5.5pt width0pt}
\def\mod#1{|#1|}
\def\modd#1{\biggl|#1\biggr|}
\def\Mod#1{\left|#1\right|}
\def\normm#1{\biggl\|#1\biggr\|}
\def\norm#1{\|#1\|}
\def\spose#1{\hbox to 0pt{#1\hss}}
\def\sub#1{^{\null}_{#1}}
\def\text #1{\hbox{\quad#1\quad}}
\def\textt#1{\hbox{\qquad#1\qquad}}
\def\T{^T\!}
\def\ntext#1{\noalign{\vskip 2pt\hbox{#1}\vskip 2pt}}
\def\enddef{{$\null$}}
\def\xint{{x_{\rm int}}}
\def\xbold{{\mathbf{x}}}
\def\zbold{{\mathbf{z}}}
\def\zbold{{\mathbf{z}}}
\def\ybold{{\mathbf{y}}}
\def\vhat{{\hat v}}
\def\fhat{{\hat f}}
\def\xhat{{\hat x}}
\def\zhat{{\hat x}}
\def\phihat{{\widehat \phi}}
\def\xihat{{\widehat \xi}}
\def\betahat{{\widehat \beta}}
\def\psihat{{\widehat \psi}}
\def\thetahat{{\hat \theta}}


\def\pthinsp{\mskip  2   mu}    
\def\pmedsp {\mskip  2.75mu}    
\def\pthiksp{\mskip  3.5 mu}    
\def\nthinsp{\mskip -2   mu}
\def\nmedsp {\mskip -2.75mu}
\def\nthiksp{\mskip -3.5 mu}


\def\Asubk{A\sub{\nthinsp k}}
\def\Bsubk{B\sub{\nthinsp k}}
\def\Fsubk{F\sub{\nmedsp k}}
\def\Gsubk{G\sub{\nthinsp k}}
\def\Hsubk{H\sub{\nthinsp k}}
\def\HsubB{H\sub{\nthinsp \scriptscriptstyle{B}}}
\def\HsubS{H\sub{\nthinsp \scriptscriptstyle{S}}}
\def\HsubBS{H\sub{\nthinsp \scriptscriptstyle{BS}}}
\def\Hz{H\sub{\nthinsp \scriptscriptstyle{Z}}}
\def\Jsubk{J\sub{\nmedsp k}}
\def\Psubk{P\sub{\nthiksp k}}
\def\Qsubk{Q\sub{\nthinsp k}}
\def\Vsubk{V\sub{\nmedsp k}}
\def\Vsubone{V\sub{\nmedsp 1}}
\def\Vsubtwo{V\sub{\nmedsp 2}}
\def\Ysubk{Y\sub{\nmedsp k}}
\def\Wsubk{W\sub{\nthiksp k}}
\def\Zsubk{Z\sub{\nthinsp k}}
\def\Zsubi{Z\sub{\nthinsp i}}


\def\subfr{_{\scriptscriptstyle{\it FR}}}
\def\subfx{_{\scriptscriptstyle{\it FX}}}
\def\fr{_{\scriptscriptstyle{\it FR}}}
\def\fx{_{\scriptscriptstyle{\it FX}}}
\def\submax{_{\max}}
\def\submin{_{\min}}
\def\subminus{_{\scriptscriptstyle -}}
\def\subplus {_{\scriptscriptstyle +}}

\def\A{_{\scriptscriptstyle A}}
\def\B{_{\scriptscriptstyle B}}
\def\BS{_{\scriptscriptstyle BS}}
\def\C{_{\scriptscriptstyle C}}
\def\D{_{\scriptscriptstyle D}}
\def\E{_{\scriptscriptstyle E}}
\def\F{_{\scriptscriptstyle F}}
\def\G{_{\scriptscriptstyle G}}
\def\H{_{\scriptscriptstyle H}}
\def\I{_{\scriptscriptstyle I}}
\def\w{_w}
\def\k{_k}
\def\kp#1{_{k+#1}}
\def\km#1{_{k-#1}}
\def\j{_j}
\def\jp#1{_{j+#1}}
\def\jm#1{_{j-#1}}
\def\K{_{\scriptscriptstyle K}}
\def\L{_{\scriptscriptstyle L}}
\def\M{_{\scriptscriptstyle M}}
\def\N{_{\scriptscriptstyle N}}
\def\subb{\B}
\def\subc{\C}
\def\subm{\M}
\def\subf{\F}
\def\subn{\N}
\def\subh{\H}
\def\subk{\K}
\def\subt{_{\scriptscriptstyle T}}
\def\subr{\R}
\def\subz{\Z}
\def\O{_{\scriptscriptstyle O}}
\def\Q{_{\scriptscriptstyle Q}}
\def\P{_{\scriptscriptstyle P}}
\def\R{_{\scriptscriptstyle R}}
\def\S{_{\scriptscriptstyle S}}
\def\U{_{\scriptscriptstyle U}}
\def\V{_{\scriptscriptstyle V}}
\def\W{_{\scriptscriptstyle W}}
\def\Y{_{\scriptscriptstyle Y}}
\def\Z{_{\scriptscriptstyle Z}}


\def\superstar{^{\raise 0.5pt\hbox{$\nthinsp *$}}}
\def\SUPERSTAR{^{\raise 0.5pt\hbox{$*$}}}
\def\shrink{\mskip -7mu}  

\def\alphastar{\alpha\superstar}
\def\lamstar  {\lambda\superstar}
\def\lamstarT {\lambda^{\raise 0.5pt\hbox{$\nthinsp *$}T}}
\def\lambdastar{\lamstar}
\def\nustar{\nu\superstar}
\def\mustar{\mu\superstar}
\def\pistar{\pi\superstar}

\def\cstar{c\superstar}
\def\dstar{d\SUPERSTAR}
\def\fstar{f\SUPERSTAR}
\def\gstar{g\superstar}
\def\hstar{h\superstar}
\def\mstar{m\superstar}
\def\pstar{p\superstar}
\def\sstar{s\superstar}
\def\ustar{u\superstar}
\def\vstar{v\superstar}
\def\wstar{w\superstar}
\def\xstar{x\superstar}
\def\ystar{y\superstar}
\def\zstar{z\superstar}

\def\Astar{A\SUPERSTAR}
\def\Bstar{B\SUPERSTAR}
\def\Cstar{C\SUPERSTAR}
\def\Fstar{F\SUPERSTAR}
\def\Gstar{G\SUPERSTAR}
\def\Hstar{H\SUPERSTAR}
\def\Jstar{J\SUPERSTAR}
\def\Kstar{K\SUPERSTAR}
\def\Mstar{M\SUPERSTAR}
\def\Ustar{U\SUPERSTAR}
\def\Wstar{W\SUPERSTAR}
\def\Xstar{X\SUPERSTAR}
\def\Zstar{Z\SUPERSTAR}


\def\psibar{\bar \psi}
\def\alphabar{\bar \alpha}
\def\alphahat{\widehat \alpha}
\def\alphatilde{\skew3\tilde \alpha}
\def\betabar{\skew{2.8}\bar\beta}
\def\betahat{\skew{2.8}\hat\beta}
\def\betatilde{\skew{2.8}\tilde\beta}
\def\ftilde{\skew{2.8}\tilde f}
\def\delbar{\bar\delta}
\def\deltilde{\skew5\tilde \delta}
\def\deltabar{\delbar}
\def\deltatilde{\deltilde}
\def\lambar{\bar\lambda}
\def\etabar{\bar\eta}
\def\gammabar{\bar\gamma}
\def\lamhat{\skew{2.8}\hat \lambda}
\def\lambdabar{\lambar}
\def\lambdahat{\lamhat}
\def\mubar{\skew3\bar \mu}
\def\muhat{\skew3\hat \mu}
\def\mutilde{\skew3\tilde\mu}
\def\nubar{\skew3\bar\nu}
\def\nuhat{\skew3\hat\nu}
\def\nutilde{\skew3\tilde\nu}
\def\Mtilde{\skew3\tilde M}
\def\omegabar{\bar\omega}
\def\pibar{\bar\pi}
\def\pihat{\skew1\widehat \pi}
\def\sigmabar{\bar\sigma}
\def\sigmahat{\hat\sigma}
\def\rhobar{\bar\rho}
\def\rhohat{\widehat\rho}
\def\taubar{\bar\tau}
\def\tautilde{\tilde\tau}
\def\tauhat{\hat\tau}
\def\thetabar{\bar\theta}
\def\xibar{\skew4\bar\xi}
\def\UBbar{\skew4\overline {UB}}
\def\LBbar{\skew4\overline {LB}}
\def\xdot{\dot{x}}

\def\Deltait{{\mit \Delta}}
\def\Deltabar{{\bar \Delta}}
\def\Gammait{{\mit \Gamma}}
\def\Lambdait{{\mit \Lambda}}
\def\Sigmait{{\mit \Sigma}}
\def\Lambarit{\skew5\bar{\mit \Lambda}}
\def\Omegait{{\mit \Omega}}
\def\Omegaitbar{\skew5\bar{\mit \Omega}}
\def\Thetait{{\mit \Theta}}
\def\Piitbar{\skew5\bar{\mit \Pi}}
\def\Piit{{\mit \Pi}}
\def\Phiit{{\mit \Phi}}
\def\Ascr{{\cal A}}
\def\Bscr{{\cal B}}
\def\Cscr{{\cal C}}
\def\Fscr{{\cal F}}
\def\Dscr{{\cal D}}
\def\Iscr{{\cal I}}
\def\Jscr{{\cal J}}
\def\Kscr{{\cal K}}
\def\Lscr{{\cal L}}
\def\Gscr{{\cal G}}
\def\Mscr{{\cal M}}
\def\Nscr{{\cal N}}
\def\lscr{\ell}
\def\lscrbar{\ellbar}
\def\Oscr{{\cal O}}
\def\Pscr{{\cal P}}
\def\Qscr{{\cal Q}}
\def\Sscr{{\cal S}}
\def\Uscr{{\cal U}}
\def\Vscr{{\cal V}}
\def\Wscr{{\cal W}}
\def\Mscr{{\cal M}}
\def\Nscr{{\cal N}}
\def\Rscr{{\cal R}}
\def\Tscr{{\cal T}}
\def\Xscr{{\cal X}}

\def\abar{\skew3\bar a}
\def\ahat{\skew2\widehat a}
\def\atilde{\skew2\widetilde a}
\def\Abar{\skew7\bar A}
\def\Ahat{\widehat A}
\def\Atilde{\widetilde A}
\def\bbar{\skew3\bar b}
\def\bhat{\skew2\widehat b}
\def\btilde{\skew2\widetilde b}
\def\Bbar{\bar B}
\def\Bhat{\widehat B}
\def\cbar{\skew5\bar c}
\def\chat{\skew3\widehat c}
\def\ctilde{\widetilde c}
\def\Cbar{\bar C}
\def\Chat{\widehat C}
\def\Ctilde{\widetilde C}
\def\dbar{\bar d}
\def\dhat{\widehat d}
\def\dtilde{\widetilde d}
\def\Dbar{\bar D}
\def\Dhat{\widehat D}
\def\Dtilde{\widetilde D}
\def\ehat{\skew3\widehat e}
\def\ebar{\bar e}
\def\Ebar{\bar E}
\def\Ehat{\widehat E}
\def\fbar{\bar f}
\def\fhat{\widehat f}
\def\ftilde{\widetilde f}
\def\Fbar{\bar F}
\def\Fhat{\widehat F}
\def\gbar{\skew{4.3}\bar g}
\def\ghat{\skew{4.3}\widehat g}
\def\gtilde{\skew{4.5}\widetilde g}
\def\Gbar{\bar G}
\def\Ghat{\widehat G}
\def\hbar{\skew{4.2}\bar h}
\def\hhat{\skew2\widehat h}
\def\htilde{\skew3\widetilde h}
\def\Hbar{\skew5\bar H}
\def\Hhat{\widehat H}
\def\Htilde{\widetilde H}
\def\Ibar{\skew5\bar I}
\def\Itilde{\widetilde I}
\def\Jbar{\skew6\bar J}
\def\Jhat{\widehat J}
\def\Jtilde{\widetilde J}
\def\kbar{\skew{4.4}\bar k}
\def\Khat{\widehat K}
\def\Kbar{\skew{4.4}\bar K}
\def\Ktilde{\widetilde K}
\def\ellbar{\bar \ell}
\def\lhat{\skew2\widehat l}
\def\lbar{\skew2\bar l}
\def\Lbar{\skew{4.3}\bar L}
\def\Lhat{\widehat L}
\def\Ltilde{\widetilde L}
\def\mbar{\skew2\bar m}
\def\mhat{\widehat m}
\def\Mbar{\skew{4.4}\bar M}
\def\Mhat{\widehat M}
\def\Mtilde{\widetilde M}
\def\Nbar{\skew{4.4}\bar N}
\def\Ntilde{\widetilde N}
\def\nbar{\skew2\bar n}
\def\pbar{\skew2\bar p}
\def\phat{\skew2\widehat p}
\def\ptilde{\skew2\widetilde p}
\def\Pbar{\skew5\bar P}
\def\phibar{\skew5\bar \phi}
\def\Phat{\widehat P}
\def\Ptilde{\skew5\widetilde P}
\def\qbar{\bar q}
\def\qhat{\skew2\widehat q}
\def\qtilde{\widetilde q}
\def\Qbar{\bar Q}
\def\Qhat{\widehat Q}
\def\Qtilde{\widetilde Q}
\def\rbar{\skew3\bar r}
\def\rhat{\skew3\widehat r}
\def\rtilde{\skew3\widetilde r}
\def\Rbar{\skew5\bar R}
\def\Rhat{\widehat R}
\def\Rtilde{\widetilde R}
\def\sbar{\bar s}
\def\shat{\widehat s}
\def\stilde{\widetilde s}
\def\Shat{\widehat S}
\def\Sbar{\skew2\bar S}
\def\tbar{\bar t}
\def\ttilde{\widetilde t}
\def\that{\widehat t}
\def\Tbar{\bar T}
\def\That{\widehat T}
\def\Ttilde{\widetilde T}
\def\ubar{\skew3\bar u}
\def\uhat{\skew3\widehat u}
\def\utilde{\skew3\widetilde u}
\def\Ubar{\skew2\bar U}
\def\Uhat{\widehat U}
\def\Utilde{\widetilde U}
\def\vbar{\skew3\bar v}
\def\vhat{\skew3\widehat v}
\def\vtilde{\skew3\widetilde v}
\def\Vbar{\skew2\bar V}
\def\Vhat{\widehat V}
\def\Vtilde{\widetilde V}
\def\UBtilde{\widetilde {UB}}
\def\Utilde{\widetilde {U}}
\def\LBtilde{\widetilde {LB}}
\def\Ltilde{\widetilde {L}}
\def\wbar{\skew3\bar w}
\def\what{\skew3\widehat w}
\def\wtilde{\skew3\widetilde w}
\def\Wbar{\skew3\bar W}
\def\What{\widehat W}
\def\Wtilde{\widetilde W}
\def\xbar{\skew{2.8}\bar x}
\def\xhat{\skew{2.8}\widehat x}
\def\xtilde{\skew3\widetilde x}
\def\Xhat{\widehat X}
\def\ybar{\skew3\bar y}
\def\yhat{\skew3\widehat y}
\def\ytilde{\skew3\widetilde y}
\def\Ybar{\skew2\bar Y}
\def\Yhat{\widehat Y}
\def\zbar{\skew{2.8}\bar z}
\def\zhat{\skew{2.8}\widehat z}
\def\ztilde{\skew{2.8}\widetilde z}
\def\Zbar{\skew5\bar Z}
\def\Zhat{\widehat Z}
\def\Ztilde{\widetilde Z}
\def\MINOS{{\small MINOS}}
\def\NPSOL{{\small NPSOL}}
\def\QPSOL{{\small QPSOL}}
\def\LUSOL{{\small LUSOL}}
\def\LSSOL{{\small LSSOL}}

\def\ds{\displaystyle}
		\def\bkE{{\rm I\kern-.17em E}}
		\def\bk1{{\rm 1\kern-.17em l}}
		\def\bkD{{\rm I\kern-.17em D}}
		\def\bkR{{\rm I\kern-.17em R}}
		\def\bkP{{\rm I\kern-.17em P}}
		\def\bkY{{\bf \kern-.17em Y}}
		\def\bkZ{{\bf \kern-.17em Z}}


		\def\beq{\begin{eqnarray}}
		\def\bc{\begin{center}}
		\def\be{\begin{enumerate}}
		\def\bi{\begin{itemize}}
		\def\bs{\begin{small}}
		\def\bS{\begin{slide}}
		\def\ec{\end{center}}
		\def\ee{\end{enumerate}}
		\def\ei{\end{itemize}}
		\def\es{\end{small}}
		\def\eS{\end{slide}}
		\def\eeq{\end{eqnarray}}
		\def\qed{\quad \vrule height7.5pt width4.17pt depth0pt} 
\def\problempos#1#2#3#4{\fbox
		 {\begin{tabular*}{0.45\textwidth}
			{@{}l@{\extracolsep{\fill}}l@{\extracolsep{6pt}}l@{\extracolsep{\fill}}c@{}}
				#1 & $\minimize{#2}$ & $#3$ & $ $ \\[5pt]
					 & $\subject\ $    & $#4$ & $ $
			\end{tabular*}}
			}
		\def\problemlarge#1#2#3#4{\fbox
		 {\begin{tabular*}{1.05\textwidth}
			{@{}l@{\extracolsep{\fill}}l@{\extracolsep{6pt}}l@{\extracolsep{\fill}}c@{}}
				#1 & $\minimize{#2}$ & $#3$ & $ $ \\[5pt]
					 & $\subject\ $    & $#4$ & $ $
			\end{tabular*}}
			}
	\def\maxproblemlarge#1#2#3#4{\fbox
		 {\begin{tabular*}{1.05\textwidth}
			{@{}l@{\extracolsep{\fill}}l@{\extracolsep{6pt}}l@{\extracolsep{\fill}}c@{}}
				#1 & $\maximize{#2}$ & $#3$ & $ $ \\[5pt]
					 & $\subject\ $    & $#4$ & $ $
			\end{tabular*}}
			}
		\def\problemsmall#1#2#3#4{\fbox
		 {\begin{tabular*}{0.55\textwidth}
			{@{}l@{\extracolsep{\fill}}l@{\extracolsep{6pt}}l@{\extracolsep{\fill}}c@{}}
				#1 & $\minimize{#2}$ & $#3$ & $ $ \\[5pt]
					 & $\subject\ $    & $#4$ & $ $
			\end{tabular*}}
			}

		\def\problem#1#2#3#4{\fbox
		 {\begin{tabular*}{0.95\textwidth}
			{@{}l@{\extracolsep{\fill}}l@{\extracolsep{6pt}}l@{\extracolsep{\fill}}c@{}}
				#1 & $\minimize{#2}$ & $#3$ & $ $ \\[5pt]
					 & $\subject\ $    & $#4$ & $ $
			\end{tabular*}}
			}
\def\maxproblem#1#2#3#4{\fbox
		 {\begin{tabular*}{0.85\textwidth}
			{@{}l@{\extracolsep{\fill}}l@{\extracolsep{6pt}}l@{\extracolsep{\fill}}c@{}}
				#1 & $\maximize{#2}$ & $#3$ & $ $ \\[5pt]
					 & $\subject\ $    & $#4$ & $ $
			\end{tabular*}}
			}
\def\maxproblemsmall#1#2#3#4{\fbox
		 {\begin{tabular*}{0.55\textwidth}
			{@{}l@{\extracolsep{\fill}}l@{\extracolsep{6pt}}l@{\extracolsep{\fill}}c@{}}
				#1 & $\maximize{#2}$ & $#3$ & $ $ \\[5pt]
					 & $\subject\ $    & $#4$ & $ $
			\end{tabular*}}
			}
\def\maxproblemsmalll#1#2#3#4{\fbox
		 {\begin{tabular*}{0.45\textwidth}
			{@{}l@{\extracolsep{\fill}}l@{\extracolsep{6pt}}l@{\extracolsep{\fill}}c@{}}
				#1 & $\maximize{#2}$ & $#3$ & $ $ \\[5pt]
					 & $\subject2\ $    & $#4$ & $ $
			\end{tabular*}}
			}

\def\largemaxproblem#1#2#3#4{\fbox
		 {\begin{tabular*}{1.0\textwidth}
			{@{}l@{\extracolsep{\fill}}l@{\extracolsep{6pt}}l@{\extracolsep{\fill}}c@{}}
				#1 & $\maximize{#2}$ & $#3$ & $ $ \\[5pt]
					 & $\subject\ $    & $#4$ & $ $
			\end{tabular*}}
			}

\def\maxxproblem#1#2#3#4{\fbox
		 {\begin{tabular*}{1.1\textwidth}
			{@{}l@{\extracolsep{\fill}}l@{\extracolsep{6pt}}l@{\extracolsep{\fill}}c@{}}
				#1 & $\maximize{#2}$ & $#3$ & $ $ \\[5pt]
					 & $\subject\ $    & $#4$ & $ $
			\end{tabular*}}
			}

\def\unconmaxproblem#1#2#3#4{\fbox
		 {\begin{tabular*}{0.95\textwidth}
			{@{}l@{\extracolsep{\fill}}l@{\extracolsep{6pt}}l@{\extracolsep{\fill}}c@{}}
				#1 & $\maximize{#2}$ & $#3$ & $ $ \\[5pt]
			\end{tabular*}}
			}

\def\unconproblem#1#2#3#4{\fbox
		 {\begin{tabular*}{0.95\textwidth}
			{@{}l@{\extracolsep{\fill}}l@{\extracolsep{6pt}}l@{\extracolsep{\fill}}c@{}}
				#1 & $\minimize{#2}$ & $#3$ & $ $ \\[5pt]
			\end{tabular*}}
			}
	\def\cpproblemsmall#1#2#3#4{\fbox
{\begin{tabular*}{0.45\textwidth}
   {@{}l@{\extracolsep{\fill}}l@{\extracolsep{6pt}}l@{\extracolsep{\fill}}c@{}}
	   #1 & & $#4 $ 
   \end{tabular*}}}

	\def\cpproblem#1#2#3#4{\fbox
		 {\begin{tabular*}{0.84\textwidth}
			{@{}l@{\extracolsep{\fill}}l@{\extracolsep{6pt}}l@{\extracolsep{\fill}}c@{}}
				#1 & & $#4 $ 
			\end{tabular*}}}
	\def\cp2problem#1#2#3#4{\fbox
		 {\begin{tabular*}{0.9\textwidth}
			{@{}l@{\extracolsep{\fill}}l@{\extracolsep{6pt}}l@{\extracolsep{\fill}}c@{}}
				#1 & & $#4 $ 
			\end{tabular*}}}
\def\z{\phantom 0}
\newcommand{\pmat}[1]{\begin{pmatrix} #1 \end{pmatrix}}
		\newcommand{\botline}{\vspace*{-1ex}\hrulefill}
		\renewcommand{\emph}[1]{\textbf{#1}}
		\newcommand{\clin}[2]{c\L(#1,#2)}
		\newcommand{\dlin}[2]{d\L(#1,#2)}
		\newcommand{\Fbox}[1]{\fbox{\quad#1\quad}}
		\newcommand{\mcol}[3]{\multicolumn{#1}{#2}{#3}}
		\newcommand{\assign}{\ensuremath{\leftarrow}}
		\newcommand{\ssub}[1]{\mbox{\scriptsize #1}}
		\def\bkE{{\rm I\kern-.17em E}}
		\def\bk1{{\rm 1\kern-.17em l}}
		\def\bkD{{\rm I\kern-.17em D}}
		\def\bkR{{\rm I\kern-.17em R}}
		\def\bkP{{\rm I\kern-.17em P}}
		\def\bkN{{\bf{N}}}
		\def\bkZ{{\bf{Z}}}

\newcommand {\beeq}[1]{\begin{equation}\label{#1}}
\newcommand {\eeeq}{\end{equation}}
\newcommand {\bea}{\begin{eqnarray}}
\newcommand {\eea}{\end{eqnarray}}
\newcommand{\mb}[1]{\mbox{\boldmath $#1$}}
\newcommand{\mbt}[1]{\mbox{\boldmath $\tilde{#1}$}}
\newcommand{\mbs}[1]{{\mbox{\boldmath \scriptsize{$#1$}}}}
\newcommand{\mbst}[1]{{\mbox{\boldmath \scriptsize{$\tilde{#1}$}}}}
\newcommand{\mbss}[1]{{\mbox{\boldmath \tiny{$#1$}}}}
\newcommand{\dpv}{\displaystyle \vspace{3pt}}
\def\texitem#1{\par\smallskip\noindent\hangindent 25pt
               \hbox to 25pt {\hss #1 ~}\ignorespaces}
\def\smskip{\par\vskip 7pt}
\def\st{\mbox{subject to}}


\title{Benefits of Multiobjective Learning in Solar Energy Prediction\\
\thanks{This Research is supported by the project AA4-11, funded by the Mathplus, Berlin, Germany (Mathplus is funded by the DFG, Deutsche Forschungsgemeinschaft).}}
\author{\IEEEauthorblockN{Aswin Kannan}
\IEEEauthorblockA{\textit{Mathematics Department} \\
\textit{Humboldt-Universität zu Berlin}\\
Berlin, Germany \\
aswin.kannan@hu-berlin.de}
}
\maketitle

\begin{abstract}
While the space of renewable energy forecasting has received significant attention in the last decade, literature has primarily focused on machine learning models that train on only one objective at a time. A host of classification (and regression) tasks in energy markets lead to highly imbalanced training data. Say, to balance reserves, it is natural for market regulators to have a choice to be more/less averse to false negatives (can lead to poor operating efficiency and costs) than to false positives (can lead to market shortfall). Besides accuracy, other metrics like algorithmic bias, RMBE (in regression problems), inferencing time, and model sparsity are also very crucial. This paper is amongst the firsts in the field of renewable energy forecasting that attempts to present a Pareto frontier of solutions (trade-offs), that answers the question on handling multiple objectives by means of using the XGBoost model (Gradient Boosted Trees). Our proposed algorithm relies on using a sequence of weighted (uniform meshes) single objective model training routines. Real world data examples from the Amherst (Massachusetts, United States) solar energy prediction panels with both triobjective (focus on accuracy) and biojective (focus on fairness/bias) classification instances are considered. Numerical experiments appear promising and clear advantages over single objective methods are seen by observing the spread and variety of solutions (model configurations).
\end{abstract}
\begin{IEEEkeywords}
pareto, XGBoost, renewable energy, forecasting, solar, fusion.
\end{IEEEkeywords}
\section{Introduction}
Multiobjective variants~\cite{multiobjco220} of the unit commitment problem are crucial in the current day world. These problems aim to consider reduction of $\textrm{CO}_2$ emissions and penetration of renewable  power in addition to cost optimization. Specifically, these take a bi-level form as follows.
\begin{align*}
\left\{P_{out}\right\}& \hspace{4mm} \min_{s \in S_w, r \in R_w} G(s,r,w) = \left[ g_1(.) \hdots g_m(.) \right].
\end{align*}
\begin{align*}
\left\{P_{in}\right\}& \hspace{4mm}\min_{\theta \in \Theta, z \in Z(\theta), w \in W(\theta,z)} F(\theta,z,w) = \left[ f_1(.) \hdots f_m(.) \right].
\end{align*}
Note that $P_{\textrm{out}}$ refers to unit commitment problem with objectives $g_i$ referring to costs, emissions etc. Here, $P_{\textrm{in}}$ denotes a machine learning problem, where the aim is to predict the amount of renewable energy that will be generated (say from solar or wind power units). We note that the problem $P_{in}$ is crucial since renewable energy has extreme variance and this helps in efficiently balancing reserves (from conventional sources of power). For the purpose of completeness, $\left\{s, r \right\}$ (from $P_{out}$) refer to generation and transmission variables and  $\left\{w, z , \theta \right\}$ (from $P_{in}$) refer to model coefficients, hyperparameters, and model/preprocessor choices respectively. We note that the focus of this paper is on the problem $P_{in}$. Here, functions $f_i$ can refer to accuracy type metrics like RMSE, MAE, $R^2$ error, F1 scores, precision, recall or computation/cost metrics like flops, model sparsity, inference time etc. We also note that other metrics pertaining to fairness/bias also form a major component of these functionals (as we will be discussing later). Taking the example of solar power, here the weights $z$ can can refer to coefficients pertaining to irradiance, humidity, and wind-speed~\cite{solarhasti21}. 

In the context of renewable sources of power, solar and wind energy have received consumer attention worldwide. 
Wind speed (an indirect quantification of the power) forecasting~\cite{classification-wind-11,classification-wind-13} has been extensively viewed and studied as a multi-label classification problem (with wind speeds falling under discrete buckets - say, 0-5 mph, 5-20 mph etc.). Multiple models, like Gradient Boosted Trees~\cite{wind-xgb21}, SVMs~\cite{svm-wind-13}, and Artificial Neural Networks~\cite{ann11-wind,ann15-wind,ann-wind-21} have been used in these cases. Regression based models have also been analyzed~\cite{regression-wind-19} in the estimation of wind speeds. While metrics like RMSE, MAE, and R-Square have been studied in this angle, it is easy to observe that these are not conflicting. Some basic works attempt to classify wind gusts (no gust vs. high gust) as means of estimating wind power production~\cite{gust11}. For a detailed survey of all such machine learning methods in lieu of wind power forecasting, the reader is asked to refer to~\cite{wind-complete-16}. In the context of solar power forecasting~\cite{solar-21-thorough}, there has been extensive work on applying regression methods with models ranging from  KNNs~\cite{knn-21-solar,huang16-knn}(K-Nearest Neighbours) and SVMs~\cite{Zulkifly2021ImprovedML,svm-alvarez-21,svm-solar-21-theo} to ANNs~\cite{ann-solar-15,ann-solar-21,germany-20-solar-ann}. Recent works using Gradient Boosted Trees have shown excellent promise~\cite{solarhasti21}. There also has been work from a classification angle~\cite{cloud-class03,cloud-class04}, where either cloud cover or irradiance levels are chosen as indirect markers of solar power. They are similarly discretized into multiple discrete buckets. For more details on readily available datasets, the reader is asked to refer to~\cite{germany-20-solar-ann,data-21-solar,usdata-18}.

A notable fact with all the above works is that they are focused on single
 objective optimization. The objectives have been some marker of accuracy (RMSE, R-Square, Balanced Accuracy, Cross Entropy Loss etc.). It can be easily noticed that the use of renewable power comes at the cost of balancing reserves (given the intermittency). Under and over forecasting are both crucial aspects that cannot be ignored. More notably, decisions have to be made quickly in realtime. This all the more necessitates the optimization of metrics pertaining to model complexity and sparsity, bias, inference time, and computational time. Some recent work has focused on looking at complexity metrics~\cite{Kommenda_2015}. But they have also been focused solely on one metric and merely consider accuracy to meet some threshold. One recent work~\cite{biobjective-17} focuses on a biobjective problem, where the interest lies in estimating the interval of generation capacity (lower and upper bounds on solar power production). The objectives of interest are the prediction probability (confidence of the interval) and the width of the interval (lower is better).
However, these works have also focused on training the model with a single standard loss and is very specific to accuracy based indicators.
The aim of this work is two-fold. Firstly, we explain the necessity of multiobjective learning in renewable power markets. Secondly, we deploy a standard decision trees (XGBoost) as our machine learning model and use weighted/scalarized optimization along the lines of~\cite{aswin2021} to train these models using custom functionals (user-defined based on the objectives). We explicitly mention that the intent of our work is to not show algorithmic benefits of weighted optimization. The aim is to merely demonstrate a working prototype of multiobjective learning for energy market predictions. This work will pave the way for future work on wind and other renewable sources, more models like neural networks, SVMs, and Random Forests, multilabel classification and regression variants of the problem, and more tailored search and first order schemes for the backend optimization (including constrained) of custom losses (and hyperparameter optimization). The rest of the paper is organized as follows. Section 2 formulates the problem and section 3 discusses about the basic algorithm in detail. Section 4 presents some promising numerical results and we conclude in section 5.
\section{Objectives and Model}
Given the penetration of renewables in the market, it is imperative that there is underlying uncertainty in generation capacities. A regulatory system operator's role is to ensure that reserves are handled appropriately such that there is no shortfall. This indirectly refers to minimization of false positives (in wind/solar power prediction). At the same time, a conservative classifier can lead to a high number of false negatives. This on the contrary can lead to very inefficient market designs with excess wastage of reserves. On a similar note, it is also important that classifers do not get biased towards a certain feature (like gender or race when it comes to algorithmic fairness). Say, a classifier predicting wind energy should not be biased towards the prediction year (say, 2017 vs 2018) or similar but different geographical locations. Noise, Model/computational complexity, and inference times are other key metrics that are crucial. The former is primarily dealt with by optimizing sparsity indicators corresponding to the model~\cite{arxiv20}. While model and computational complexity are heavily tied to the hyperparameters of the model, inference times are blackbox based and require special attention. Even when the goal is to maximize accuracy, the choice of the related accuracy metric is very much tied to the setting of interest. Say, smaller planning horizons can demand higher importance to metrics like balanced accuracy, while longer horizons can require optimization of AUCROC or RMBE (Root Mean Bias Error). For a detailed choice of metrics, please refer to~\cite{solar-21-thorough}. Our multiobjective machine learning problem can be compactly stated as follows.
\begin{align*}
& \min_{z,w}  \hspace{4mm} F(y^p(z,w,x^{l}),y^{l}) = \pmat{f_1(y^p(z,w,x^{l}),y^{l}) \\ \vdots \\ f_m(y(z,w,x^{p}),y^{p})}^T \\
& \text{Subject To:} z \in Z, \hspace{2mm} w \in W(z).
\end{align*}
Here, $f_i(.)$ can represent any objective like balanced accuracy, precision, recall, RMSE, R-Square, AUCROC, complexity, RMBE, or inference time etc. depending on the problem type. Standard losses like cross entropy or hinge can also be used.  Note that $y^p$ and $y^l$ denote the respective model predictions and training data labels. Here, $x^l$ refers to the feature vector in the training data. Most importantly, $z$ refers to model hyperparameters and $w$ refers to the weights or model parameters. Say, if Neural Networks are deployed, z corresponds to the number of hidden layers and neurons, while w refers to the weights associated with each of these layers/neurons. In cases of decision trees, z can refer to tree depth and the number of leaves. The machine learning problem can include the case of data pre-processors and feature engineering components. This framework is referred to as pipeline optimization and our model mentioned above will follow without any loss of generality. However, pipeline optimization is not in the scope of our current research. We focus on three classes of problems (as to be discussed in the numerics section): triobjective and biobjective problems with focus on accuracy and thirdly a biobjective problem with focus on bias and accuracy. It can also be noted that forecasting problems can further be complicated by means of constraints on predictions, metrics, and model parameters~\cite{LUO2021-constraints-solar}. However, we note that this is also not in the scope of our current work.

\section{Algorithms}
\begin{figure}[t!]
\centering
\includegraphics[width=0.8\linewidth]{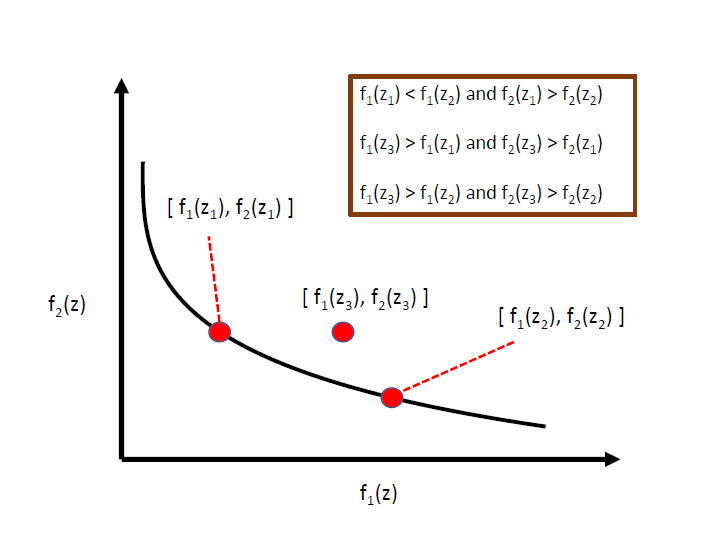}
\caption{Sample Pareto frontier for a biobjective optimization problem.}
\label{fig:pareto}
\end{figure}

Multiobjective learning has recently received good attention in literature, particularly in multi-task learning~\cite{ruder2017overview,zhang2017survey,Sener2018MultiTaskLA}. Note that the goal with multi-task learning is to simultaneously learn model parameters $w$ that aim at performance improvement. However most of these problems are triggered by performance measures over multiple datasets and a single objective. Recent advancements have handled multiobjective learning problems by means of generating a Pareto frontier through black-box optimization methods. These can be vaguely classified under two subsets, namely Genetic algorithms and Bayesian optimization. Genetic algorithms~\cite{nsga2,genetic2011inria} are population based and use Pareto dominance as their criteria to update solutions/iterates. Multi-objective Bayesian Optimization (MOBO) methods~\cite{usemo,konakovic2020diversity,daulton2020differentiable} build surrogate models for all objectives of interest and select iterates in batches based on the expected improvement of the hypervolume of the Pareto frontier (also an estimate from the sequential information on iterates/objectives). The hypervolume indicator~\cite{audet18,stk12} is a standard metric that has been well accepted in literature for pareto based multiobjective optimization. For clarity and completeness, we define a Pareto frontier as per figure~\ref{fig:pareto}. A simple biobjective problem is shown in this case. It can easily be seen that points $x_1$ and $x_2$ dominate point $x_3$. However, it can also be observed that moving from point $x_1$ to $x_2$, does not guarantee a descent in both objectives $f_1$ and $f_2$. These are referred to non-dominated points and the set of such points constitute the Pareto frontier.
We note that all these approaches are however generic and assume that there is little information on the structure of the objectives. However in the case of machine learning, we note that the structure of the underlying model is well defined. Additionally, the gradients and hessians of the objectives in terms of the model predictions are mostly known and given by closed form expressions. It is hence obvious to note that the dependency of the objectives on the model parameters (weights) is not completely blackbox based and more information can be extracted from the above knowledge (efficient exploitation of structure). We also note that most of these algorithms also focus on hyperparameter optimization (in the space of $z$), where multiple configurations of $z$ define the Pareto frontier. However, for each configuration in $z$, the models are trained only with standard losses (like cross entropy) which are only accuracy oriented. In short, there is very little in terms of carving out a frontier in terms of multiple configurations of $w$'s for a fixed configuration of $z$. Work along the direction of search based methods~\cite{search2011, Fejes16} have also proved to be very efficient, but they are also restricted to hyperparameter optimization. Recent works exploit this aspect and train on custom losses~\cite{aswin2021}. These custom losses are constructed based (MPO-Model Parameter Optimization) on sequential scalarization of the objectives by means of weights ($\rho$). These weights are chosen in an adaptive sense, based on the evolving frontier of optimal solutions (in the objective space). In addition, Bayesian solvers like USEMO are deployed for the hyperparameter tuning (HPO-Hyperparameter Optimization) and the final frontier is a result of the fusion of the HPO and MPO schemes.

In this work, we use similar Bayesian methods for the HPO portion. However, for the MPO portion, we resort to a standard uniform weights based approach. Say, with two objectives, the following set of weights presents a uniform choice: $\left\{ \rho = (0,1), (0.25,0.75), (0.5,0.5),(0.75,0.25), (1,0). \right\}$.
The optimization problem of our interest is a bi-level case stated as follows.
\label{prob:prob1}
\begin{align}
\underbrace{\min_{z}}_{\text{HPO}} & \left\{ \underbrace{\min_{w \in \mathcal{W}(z)}}_{\text{MPO}} \quad \sum_{j=1}^m \rho_j f_j(w) \right\} \\
& \text{ subject to } z \in Z.
\end{align}
Here, note that $W(z)$ refers to the set of possible model parameters that can be chosen, given a hyperparameter configuration $z$. As mentioned earlier, $z$ corresponds to the number of layers and neurons in the case of neural networks. It can refer to even the choice of activation functions and pooling layers. In the case of XGBoost, this refers to the maximum number of leaves, depth of a tree, and number of trees. This can also include algorithmic hyperparameters like learning rate or classification threshold. It is easy to observe that the space of $z$ can be marked by categorical/integer variables. The weights $w$ usually are not constrained by themselves.
Say, if there are 4 hidden layers with 10 neurons each, we would have 40 weights, (i.e) $|W(z)| = 40$ (cardinality denoted by |.| in this case). With objectives on model complexity, sparsity becomes very crucial. In such cases, one of the objectives takes the form, say $f_j = \| w\|_2$. For completeness, we re-define the parametrization of the objectives. Say, let $f_1$ denote MSE (Mean Square Error). Then,
$$f_j(.) = \sum_{i=1}^{N}\|y^{p}(z,w(z),x^i)-y^i \|^2 .$$
Other metrics and indicators are defined accordingly. We summarize our algorithm as follows for the purpose of clarity.
\begin{itemize}
\item Step 1:  Solve the MOO as-is by means of Bayesian optimization that treats the objectives as simulation based (derivative-free/blackbox) for a fixed time budget ($t_0$). Let $\left\{ z_k \right\}_{k=0}^{K}$ denote the resulting trajectory of hyperparameters.
\item Step 2: Fix the accuracy threshold $\theta$. For all $k$ such that $f_{acc}(z_k) \geq \theta$, sort the set of iterates w.r.t time (computational time).
\item Step 3: Till time budget $t-t_0$ is attained, perform MPO (inner problem) by fixing the configurations in the sorted order. For each configuration, generate a set of weights from the uniform space as denoted earlier. Another example is given below.
\end{itemize}
Say, for a triobjective problem with 10 mesh points allowed for evaluation, the weights $\rho$ are $\left\{ 0, 1/3, 2/3\right\}$,  $\left\{ 0, 2/3, 1/3\right\}$. $\left\{ 1/3, 2/3, 0\right\}$, $\left\{ 2/3, 1/3, 0\right\}$, $\left\{ 1/3, 0, 2/3\right\}$, $\left\{ 2/3, 0, 1/3\right\}$, $\left\{ 1/3, 1/3, 1/3\right\}$, $\left\{ 1, 0, 0\right\}$, $\left\{ 0, 1, 0\right\}$, and $\left\{ 0, 0, 1\right\}$. This refers to a total $10 = 3^2+1$ number of weights. That is, with mesh size $n = 3$, the total number of weights $= n^2+1$.

\paragraph{Non-weighted methods:} It can be noted that weighted/scalarization methods do not show a similar promise when the underlying pareto frontier happens to be nonconvex. In such settings, $\epsilon$ constraint based methods~\cite{epsilon71} have demonstrated good computational performance. These rely on solving a sequence of single objective optimization problems, say as follows (triobjective case shown below).
\begin{align*}
\min_{w} & \hspace{3mm} f_1(w) \\
\st & \hspace{3mm} f_2(w) \leq \epsilon_2 \\
& \hspace{3mm} f_3(w) \leq \epsilon_3, \hspace{3mm} w \in W(z).
\end{align*}
This is recursively done over several careful choices of $\epsilon$ and objectives. Besides, search based methods like Nelder-Mead~\cite{Fejes16} and Multi-Mads~\cite{mads21multiobj} have been deployed in the context of HPO without using weights. These can be extended to MPO if some knowledge on the expected solutions is known based on subject matter expertise (else, they can turn out to be expensive). Hybrid regimes have also been proposed more recently, which use gradients and search techniques together~\cite{AugLag14}. Lastly, reference point based methods that use Shapley values (Explainable AI) have shown promise in estimating approximate pareto frontiers~\cite{shapley22moo}, when problem dimensions are not huge. These non-weighted methods prove to be one useful future research direction that has an immense potential for improving the quality of solutions/Pareto frontier.
\section{Numerics}
\begin{figure*}[t!]
\centering
\begin{subfigure}{.44\textwidth}
  \centering
  \includegraphics[height=0.55\linewidth]{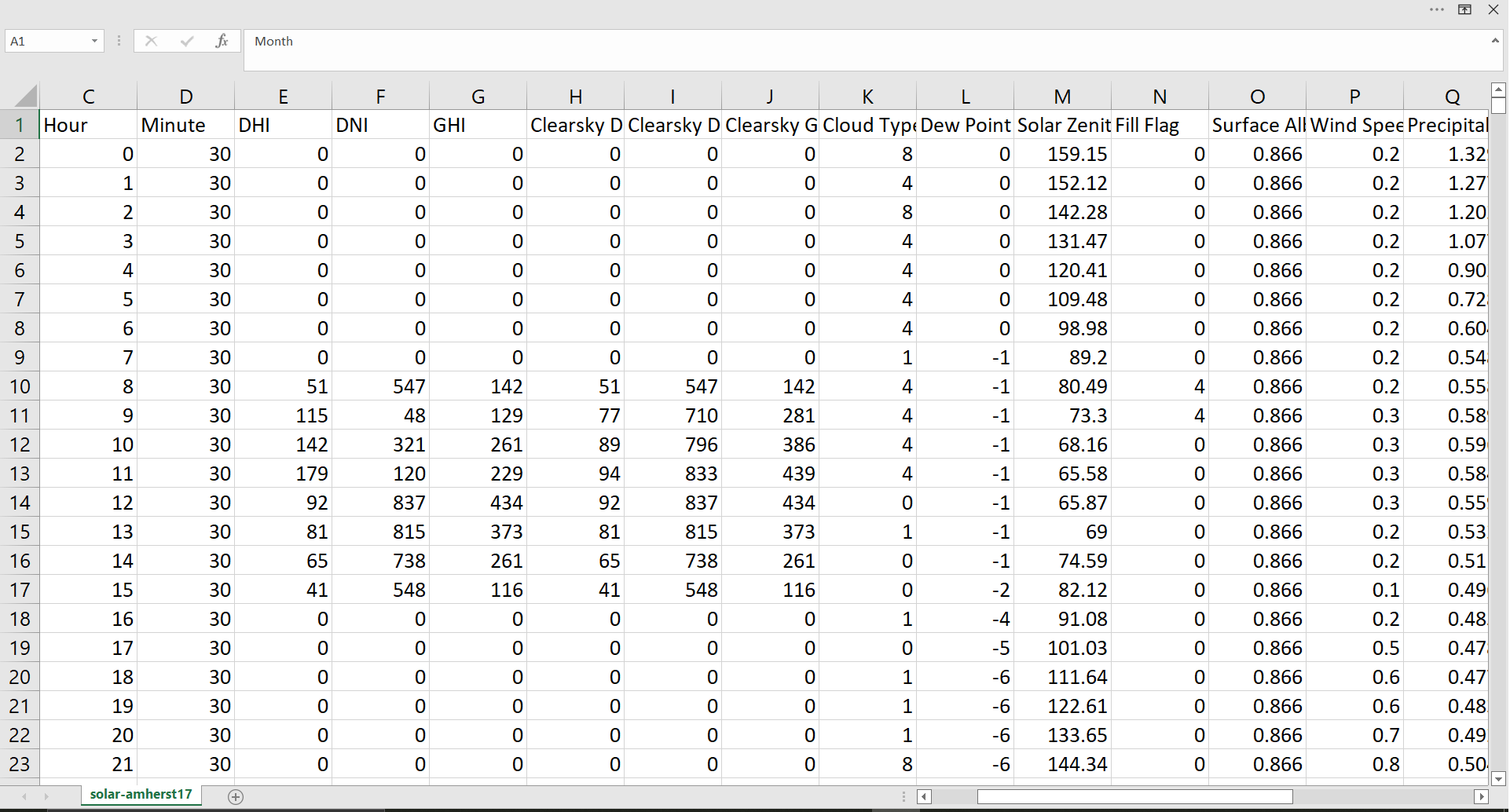}
\end{subfigure}%
\hspace{3mm}
\begin{subfigure}{.44\textwidth}
  \centering
  \includegraphics[height=0.55\linewidth]{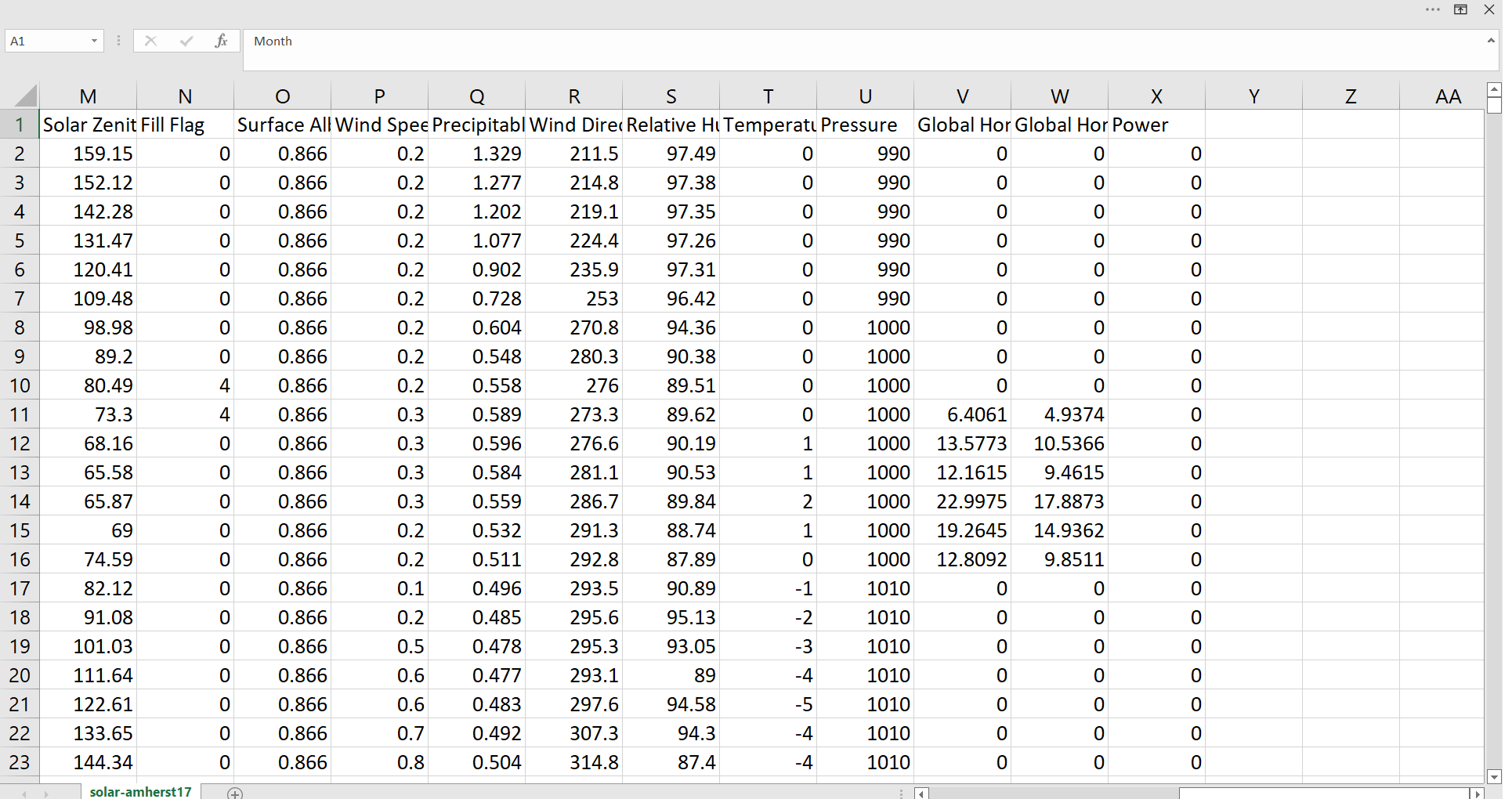}
\end{subfigure}
\caption{Visualization of Data.}
\label{fig:data}
\end{figure*}

We present our results with the University of Massachusetts dataset (UMass)~\cite{solarhasti21,umass-cs-solar,nrel-data}. This presents real-time observations of solar incidence features and generations at the Computer Science Department at UMass. We specifically note that this is a two-step process, where data on the features like incidence angle, GHI, wind-speeds, pressure, and temperature are obtained from ``weather.com'' and the data pertaining to the actual power production is fetched from the University of Massachusetts. Also, it is to be noted that UMass houses six generation facilities (for PV-power production) and we just consider the one from the Computer Science department for our work. We also note that some data points are irrelevant and we drop them during our pre-processing stage (power generation is reported only from certain hours of the day). Additionally, we note that we use a ``toy'' version of the dataset, where the labels are converted/encoded to be suitable for classification analysis. That is, positive power generation is labeled as ``1'' and no power generation is labeled as ``0''. For ease of parsing and understanding, a snapshot of the finally used data is shown in figure~\ref{fig:data}. While regression models can also be studied (given the inherently continuous data labels), we leave this for future research. This is mainly because, most of the renewable energy literature has focused on either binary or multi-label classification.
We present two sets of results, one on using different metrics of accuracy and one tied to bias/fairness. Each of these sets further expand onto two different problems each. All the computation was done on a cluster with 12 cores and processor speed \& memory of 3.1 Ghz \& 32 GB respectively. All our results are presented with the XGBoost model.
\subsection{Accuracy Indicators}
\begin{figure}[t!]
\centering
\begin{subfigure}{.23\textwidth}
  \centering
  \includegraphics[height=0.83\linewidth]{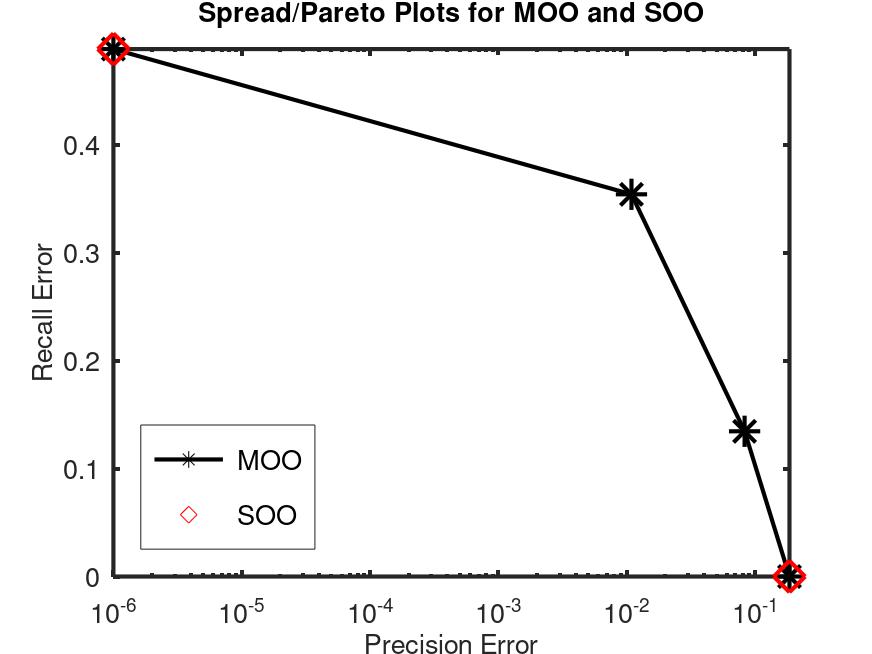}
  \caption{Pareto frontier.}
  \label{fig:fusionsort}
\end{subfigure}%
\begin{subfigure}{.23\textwidth}
  \centering
  \includegraphics[height=0.83\linewidth]{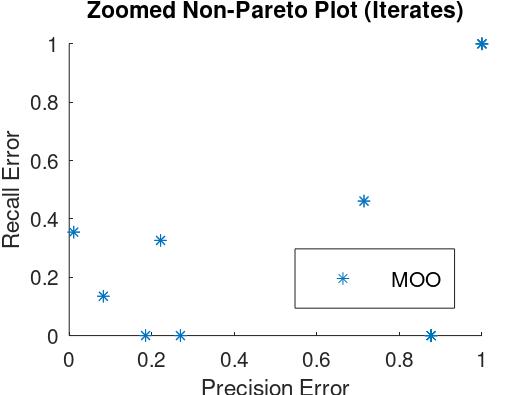}
  \caption{Scatter - Zoomed.}
  \label{fig:triobj}
\end{subfigure}
\caption{Performance of Single Objective (SOO) and Multiobjective Optimization (MOO) approaches on the biobjective problem.}
\label{fig:biobj}
\end{figure}

Two period power markets with increased solar power penetration are characterized by increased commitments and participation (across all firms) in the forward market ~\cite{kannanoms13}. Deviations from forward market bids naturally face huge penalties from market operators. Given the underlying uncertainty in solar power, this further necessitates crucial planning of the reserve market. When planning with reserves, mere overall accuracy of the solar power predictions do not suffice. Given the highly imbalanced nature of the datasets (intermittency of power generation), both false positives and negatives require minimization. False positives can lead to a higher risk in shortfall and outages. At the same time, false negatives can also increased purchases of reserves and significantly higher costs. Outages are usually balanced by emergency supplies, which come at a significantly higher cost and international trade relations. The balance between these costs differ from one country to another, year to year, and even season to season. While joint minimization of both false negatives and positives are ideal, such a setting cannot be achieved and we have to resort to Pareto solutions. Practitioners can pick solutions (models) that are usable according to their dynamic needs. For this study, we use a concise version of the data from both 2017 and 2018 (with nearly 6000 records - after removing some irrelevant records). In the given dataset, we note that 9-10 percent of the data presents with positives (power generation) and nearly 90 percent of the time, there is no power generation. We consider two years worth of data (2017 and 2018). Here, some of the features represent DHI, DNI, GHI, Clearsky parameters, Cloud cover, Dew Point, wind speed, temperature, and pressure. For a more detailed note on the features, the reader is asked to refer to~\cite{solarhasti21}, which uses the same dataset (in a slightly different context). Here, ~\cite{solarhasti21} presents a regression based single-objective approach to predict accurate models for solar power forecasting (also using XGBoost). In our case, we re-label the actual generation data to binary levels (which is very common in solar and wind energy markets). Say, generation levels of 0 and 1.5 MW are modified to 0 and 1 respectively (No/Yes labels). As much as we notice in regression problems, no classifier is free from bias. In this process, we also consider mean bias error as one of our objectives. We consider the following seven  hyperparameters for our study. The associated lower and upper bounds are also enclosed.
\begin{itemize}
\item Architectural (integer): tree depth (1 to 13), maximum leaves (1 to 100), and number of rounds (1 to 100).
\item Algorithmic/fitting hyperparameters (continuous): classification threshold (0 to 1) and learning rate (0 to 1).
\item Others (continuous): minimum child weight (0 to 10), max-delta-step (0 to 10).
\end{itemize}
\paragraph{Precision and Recall}
First, we consider the case of just false positives and negatives. In this work, we primarily focus on model parameter optimization. While the importance of hyperparameter optimization is truly noted, we generate just three such configurations from four equally spaced intervals (say, L+U/4, L+U/2, and 3L+U/4 if lower and upper bounds are represented by L and U). The objective of this work is to compare the pareto frontier (obtained by using weighted optimization) with single objective solutions (obtained by training with cross entropy losses). We note that multiple losses and hyperparameter (including the three above, but more) configurations are considered for the case of single objective optimization. For the MPO portion, we generate weights from a uniform distribution, 5 equally spaced points in the 2D space, between 0 and 1. That is, $w_1 = \frac{i}{4}$ and $w_2 = 1-w_1$, for $i = 0, 1, 2, 3$, and $4$. Figure~\ref{fig:biobj} shows the s spread plots for both our weighted scheme and a single objective examination based approach. We note that the number of evaluations (number of times models are trained) in each case are 15 (for consistency purposes).
We also note that the actual expense of evaluating the function is higher in the case of using HPO. We report results from both years 2017 and 2018 for the purpose of completeness. For completeness, our objectives are defined as follows.
\begin{align*}
&{\tt prec}(\mathbf{y^a},\mathbf{y^p}) = \frac{\sum_{i=1}^m (1-y_i^a)y_i^p}{\sum_{i=1}^{m}y_i^p}, \\
&{\tt rec}(\mathbf{y^a},\mathbf{y^p}) =  \frac{\sum_{i=1}^m  (1-y_i^p)y_i^a}{\sum_{i=1}^{m}y_i^a}.
\end{align*}
\paragraph{Triobjective Problems}
In this case, we add another accuracy metric, namely the Mean Bias Error and consider the triobjective version of the problem. For completeness, we define MBE as follows.
$$MBE = \frac{1}{N}\left|\sum_{i=1}^{N} \left\{ y_i^p - y_i^a \right\}\right|.$$
Note that $y_i^p$ and $y_i^a$ refer to the model predictions and the actual labels respectively. For this triobjective problem, we just show the performance comparison of both approaches with increasing function evaluations. For the mesh generation with weights, we assume the following pattern, $w_1 = \frac{i}{3}, w_2 = \frac{j}{3}$, and $w_3 = 1-w_1-w_2$ (for $i = 0, 1, 2$, and $3$ and $j$ defined similarly). Figure~\ref{fig:triobj} shows the comparison of the hypervolumes across the single and multiobjective approaches. Note that single objective optimization is performed using a standard cross entropy loss (other losses or metrics like precision/recall can also be used). The benefits of using multiobjective optimization cannot be very easily noticed in this case. The bias metric does not seem to be conflicting and it can be inferred that there is very little bias in the dataset/classifier combination.
\begin{figure}[t!]
\centering
\includegraphics[height=0.75\linewidth]{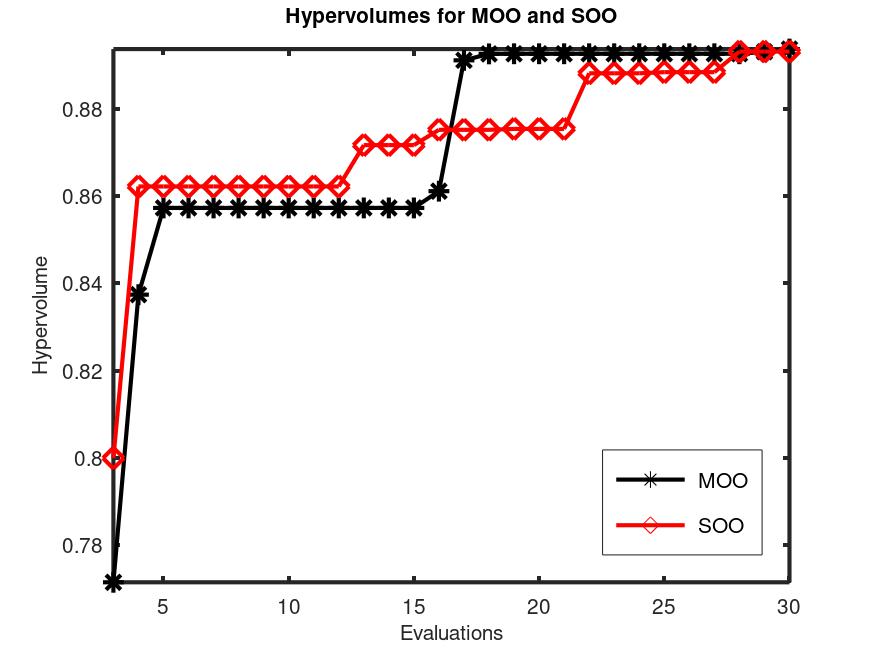}
\caption{Triobjective Problem - comparison of MOO/SOO.}
\label{fig:triobj}
\end{figure}
\subsection{Fairness Indicators}
While accurate classifiers are good, research over the last decade has shown that these can be unduly biased towards some features, to enforce accuracy. Say, classifers can tend to give undue importance to race or gender to get accurate predictions. Fairness in AI has been very deeply studied in the last decade~\cite{gajane,calders,aif360_datasets}. In the context of solar power generation, it can easily be noticed that some parameters do not contribute much towards predictions or classification tasks. Say, generation levels are not bound to change significantly within years (say 2017and 2018) if other features remain the same. Also in our dataset, features like wind speeds are essential, but not sole markers of production levels. In such settings, there is always a necessity for practitioners to compromise slightly in accuracy by ensuring that the classifiers do not get biased significantly. The extent of balancing these objectives is highly subjective and therefore leads to a similar biobjective optimization problem. For our fairness related problems, we use the Calders-Verwer or the CV-score (also related to DPD - Disparate Parity Difference) measure as our indicator.
The DPD measure represents the gap between the probabilities of getting positive outcomes by a binary predictor $\hat{Y}$ for two different sensitive groups in the input data. It is defined as:
\[DPD := \big \lvert Pr(\hat{Y}=+|A=1) - Pr(\hat{Y}=+|A=0) \big \rvert \]
Here, \emph{A} denotes a binary sensitive attribute or the fairness feature (e.g., wind-speed or year) in the input data with A=1 and A=0 representing privileged and unprivileged groups respectively. Predicted outcome $\hat{Y}$ can be positive (favorable) or negative (unfavorable). A fair predictor will have similar proportions of positive outcomes within each sensitive group resulting in a lower DPD value. Noting that $A$ marks the binary sensitive attribute, we define index sets $I_A$ and $J_A$ referring to the sample points (feature data) that correspond to 0's and 1's of the attribute. The objective takes the following form:
\begin{align*}
&I_A = \left\{j| m \geq j \in Z,  A_j = 1 \right\}, \\
&J_A = \left\{j| m \geq j \in Z, A_j = 0 \right\} \\
&{\tt DPD}(\mathbf{y^a},\mathbf{y^p}) = \left| \frac{\sum_{i \in I_A} y_i^p}{|I_A|_c} - \frac{\sum_{i \in J_A} y_i^p}{|J_A|_c} \right|.
\end{align*}
Note that the function $|.|_c$ represents the cardinality of a set, whereas $|.|$ denotes the absolute value function. Finally, the Balanced Accuracy Error metric (which is 1-Balanced Accuracy) is quite standard and is defined as:
$$ {\tt bal\_err}(\mathbf{y^a},\mathbf{y^p}) = \frac{1}{2|P|_c}\sum_{i \in P} |1-y_i^p | + \frac{1}{2|N|_c}\sum_{i \in N} |y_i^p |.$$
Also, $P$ and $N$ denote the set of points (labels) classified as positive (1's) ad negative (0's) respectively.
\paragraph{Wind Speed}
Say, as an example the user/subject matter expert knows that changes in wind speeds by 20 mph does not cause any change in the amount of energy trapped by the solar panel (given other features remain the same). If a classifier predicts elsewise, it is deemed to be biased. In this exercise, we analyze the fairness of the classifier by fixing the wind-speed column.
Figure~\ref{fig:fairwind} shows the spread plots for both single and multiobjective optimization. Some improvement and the variety of solutions can be  noticed when using multiobjective optimization. The number of intermediate points in the frontier will improve with increasing more weights. We note that the purpose of this exercise is to just compare the difference between single and multiobjective optimization. We do not intend to show a performance improvement in using one scheme over the other. Our choice of weights and hyperparameters is the same as earlier. The same exercise can be repeated with other less relevant features to fine tune our process of computing efficient models.
\begin{figure}[t!]
\centering
\begin{subfigure}{.23\textwidth}
  \centering
  \includegraphics[height=0.83\linewidth]{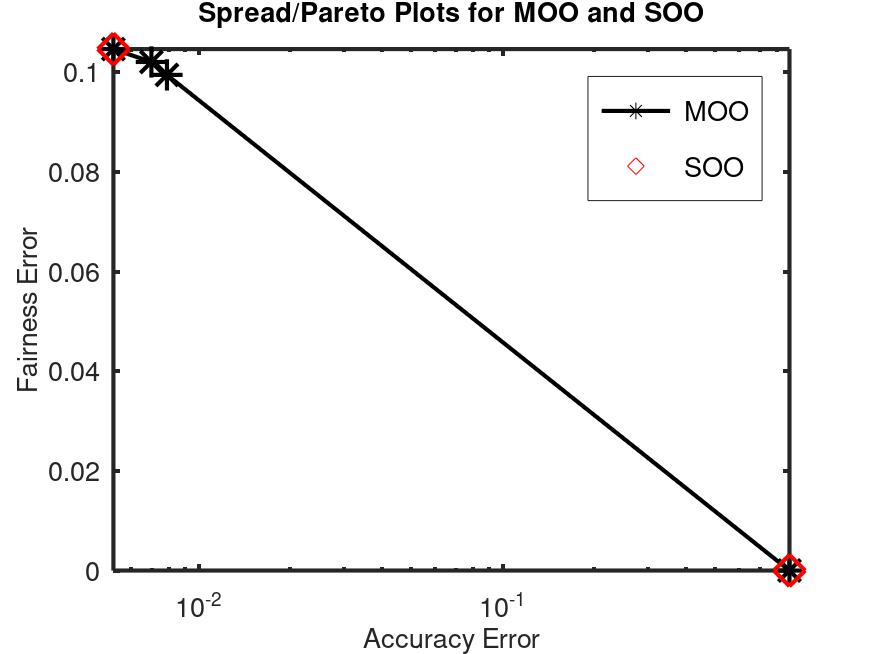}
  \caption{Pareto frontier.}
  \label{fig:fusionsort}
\end{subfigure}%
\begin{subfigure}{.23\textwidth}
  \centering
  \includegraphics[height=0.83\linewidth]{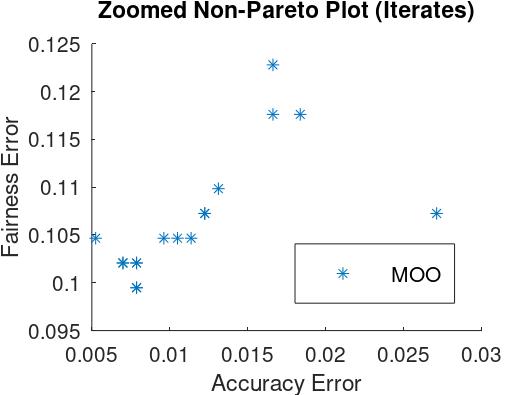}
  \caption{Scatter - Zoomed.}
  \label{fig:fairwind}
\end{subfigure}
\caption{Performance of Single Objective (SOO) and Multiobjective Optimization (MOO) approaches on the Bias/Fairness problem with indicator ``Wind-Speed''.}
\label{fig:fairwind}
\end{figure}

\begin{figure}[t!]
\centering
\begin{subfigure}{.23\textwidth}
  \centering
  \includegraphics[height=0.83\linewidth]{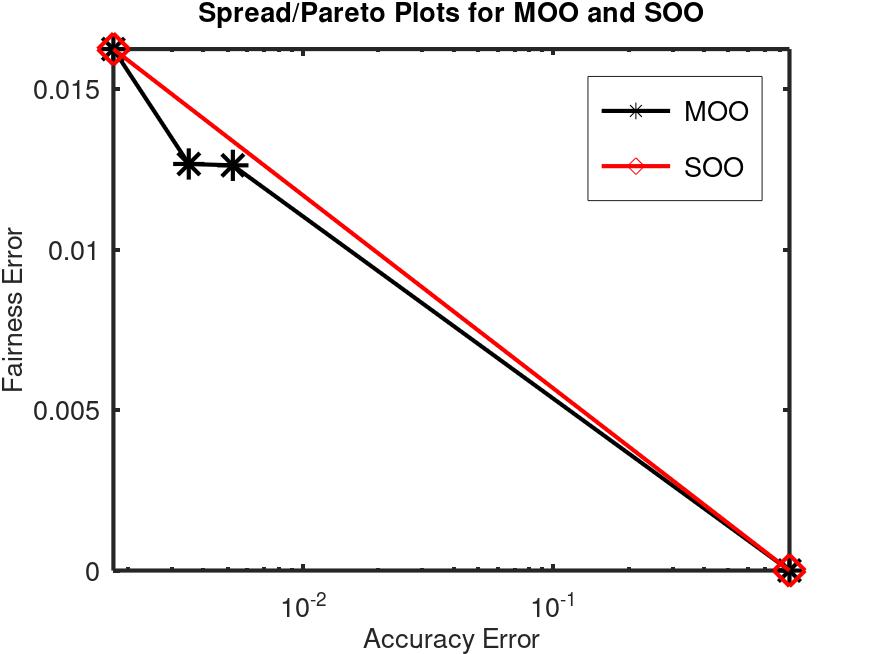}
  \caption{Pareto frontier.}
  \label{fig:fusionsort}
\end{subfigure}%
\begin{subfigure}{.23\textwidth}
  \centering
  \includegraphics[height=0.83\linewidth]{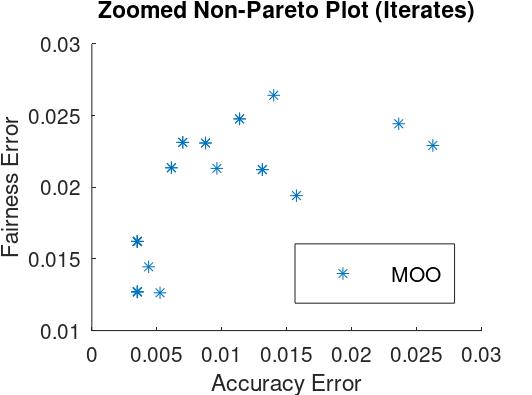}
  \caption{Scatter - Zoomed.}
  \label{fig:fairyear}
\end{subfigure}
\caption{Performance of Single Objective (SOO) and Multiobjective Optimization (MOO) approaches on the Bias/Fairness problem with indicator ``Year''.}
\label{fig:fairwind}
\end{figure}

\paragraph{Year}
In this exercise, we additionally add the "year" column as one of the features. As a hypothesis, it can be stated that under usual circumstances, the solar energy outputs between any two years will be similar. In other words, classifiers can compromise on  this hypothesis in an attempt to yield accurate models. We hence repeat the exercise with DPD with year as the fairness feature. Similarly, figure~\ref{fig:fairyear} shows a comparison in spread and hypervolume between single and multiobjective solutions. The message conveyed remains the same. This all the more necessitates the essence of training custom losses, with user defined metrics.
\section{Conclusion}
The paper presented a basic technique to present Pareto solutions for solar energy forecasting problems. Benefits of training on multiple objectives (in comparison to single objective) can easily be observed from the numerical results presented. It can be noted that the presented approach is tied to a stylized case with the XGBoost model and optimization dependent on scalarization techniques (model parameter training). We note that this work forms an application base that strongly cements the need for for comprehensive future research along two different directions. The first direction refers to the use of more models like MLP, Neural Nets, SVMs, and Random Forests, extensions to wind power problems and regression based settings, and use of more datasets (say from Open Weather Map). The second direction should focus on the deployment of Bayesian methods and Direct Search type search methods (say, MADS or Nelder-Mead) for hyperparameter tuning and Frank-Wolfe type methods for addressing other complex metrics like model inference-time and sparsity.
\section*{Acknowledgment}
The author would like to thank the Mathplus (funded by the DFG, Deutsche Forschungsgemeinschaft) for the project grant AA4-11, titled "Using Mathematical Programming to Enhance Multiobjective Learning". The author would also like to thank anonymous reviewers for useful constructive comments on the application and algorithms.
\bibliographystyle{IEEEtran}
\bibliography{sample}
\end{document}